%Basic AmsTeX article template

\input amstex

\define\caf{$\Cal F$}
\define\cafi{${\Cal F}_i$}
\define\tcaf{$\widetilde{\Cal F}$}
\define\tim{$\widetilde{M}$}
\define\tix{$\widetilde{X}$}
\define\gam{$\gamma$}
\define\lam{$\lambda$}
\define\tgam{$\widetilde{\gamma}$}
\define\xo{$x_0$}
\define\txo{$\widetilde{x_0}$}

\define\del{$\partial$}

\define\cafu{${\Cal F}_1$}
\define\caft{${\Cal F}_2$}
\define\cl{$\Cal L$}

\define\gampr{$\gamma^\prime$}
\define\ctln{\centerline}
\define\ssk{\smallskip}
\define\msk{\medskip}

\overfullrule=0pt

\documentstyle{amsppt}

\topmatter

\title\nofrills Tautly foliated 3-manifolds with no 
$\bold R$-covered foliations \endtitle
\author  Mark Brittenham \endauthor

\leftheadtext\nofrills{Mark Brittenham}
\rightheadtext\nofrills{No $\bold R$-covered foliations}

\affil   University of Nebraska \endaffil
\address   Department of Mathematics and Statistics, University 
of Nebraska, Lincoln, NE 68588-0323\endaddress
%\curraddr           \endcurraddr
\email   mbritten\@math.unl.edu \endemail
%\dedicatory       \enddedicatory
%\date                    \enddate 
\thanks   Research supported in part by NSF grant 
\# DMS$-$9704811 \endthanks
%\translator         \endtranslator
%\subjclass            \endsubjclass   

\abstract We provide a family of 
examples of graph manifolds which
admit taut foliations, but no 
$\bold R$-covered foliations.
We also show that, with very few 
exceptions, $\bold R$-covered foliations
are taut.\endabstract

\endtopmatter 

\heading{\S 0 \\ Introduction}\endheading

A foliation \caf\ of a closed 3-manifold 
$M$ can be lifted to a foliation \tcaf\
of the universal cover \tim\ of $M$; if 
the foliation \caf\ has no Reeb components,
the leaves of this lifted foliation are 
all planes, and Palmeira [Pa] has shown that \tim\
is homeomorphic to $R^3$. Palmeira also 
showed that \tcaf\ is homeomorphic to
(a foliation of $R^2$ by lines)$\times R$. 
The space of leaves of \tcaf, the quotient
space obtained by crushing each leaf to a
point, is homeomorphic to the space of
leaves of the foliation of $R^2$, and is a 
(typically non-Hausdoff) simply-connected
1-manifold. 

If the space of leaves is Hausdorff (and 
therefore homeomorphic to $\bold R$), we say
that the foliation \caf\ is {\it 
$\bold R$-covered}. Examples of $\bold R$-covered foliations
abound, starting with surface bundles 
over the circle; the foliation by fibers is 
$\bold R$-covered. Thurston's notion 
of `slitherings' [Th] also provide a large 
collection of examples. A great deal has 
been learned in recent years about 
$\bold R$-covered foliations and the manifolds 
that support them (see, e.g., [Ba],[Ca],[Fe]),
especially in the case when the underlying 
manifold $M$ contains no incompressible
tori.

The purpose of this note is to provide 
examples of 3-manifolds which admit taut foliations, 
but which do not admit any $\bold R$-covered 
foliations. All of our
examples are drawn from graph manifolds, 
and so all contain an incompressible torus,
and, in fact, a separating one. The technology 
does not exist at present to identify
examples which are atoroidal, and it is 
perhaps not clear that such examples should
be expected to exist. That we can find 
the examples we seek among 
graph manifolds relies on two facts: (a) we have a good
understanding [BR] of how taut and 
Reebless foliations can meet an incompressible
torus, and (b) we understand [Br1],[Cl],[Na] 
which Seifert-fibered spaces can contain
taut or Reebless foliations. These two 
facts have been used previously [BNR] to find graph
manifolds which admit foliations with 
various properties, but no {\it ``stronger''}
properties; the examples we provide here 
are in fact the exact same examples used
in [BNR] to illustrate its results. We will 
simply look at them from a
somewhat different perspective.

This paper can, in fact, be thought of as a 
further illustration of this `you can get
this much, but no more' point of view towards 
foliating manifolds. Namely, you can
get taut, but not $\bold R$-covered. 
$\bold R$-covered is, it turns out, almost (but not
quite) always a stronger condition: in all but 
a very small handful of instances, an $\bold R$-covered
foliation must be taut, as we show in the next 
section. This result seems to have been
implicit in much of the literature on $\bold R$-covered 
foliations; a different proof of this 
result, along somewhat different lines, can be found in [GS].

\heading{\S 1 \\ $\bold R$-covered almost implies taut}\endheading

In this section we show that in all but a very 
few instances, an $\bold R$-covered foliation
must be taut. We divide the proof into two parts; 
first we show that a Reebless $\bold R$-covered
foliation is taut, and then describe the manifolds 
that admit $\bold R$-covered foliations with Reeb components.

Recall that a {\it Reeb component} is a solid torus 
whose interior is foliated by planes transverse
to the core of the solid torus, each leaf limiting 
on the boundary torus, which is also a leaf.
(There is a non-orientable version of a Reeb 
component, foliating a solid Klein bottle, which we
will largely ignore in this discussion. It can 
be dealt with by taking a suitable
double cover of our 3-manifold.)
We follow standard practice and refer to both the solid 
torus and its foliation as a Reeb component.
A foliation that has no Reeb components is called 
{\it Reebless}. A foliation is {\it taut} if for
every leaf there is a loop transverse to the 
foliation which passes through that leaf.
Taut foliations are Reebless.

\proclaim{Lemma 1} If a closed, irreducible 
3-manifold $M$ admits a Reebless,
$\bold R$-covered foliation \caf\ containing a
compact leaf $F$, then every component of $M|F$, the
manifold obtained by splitting $M$ open
along $F$, is an $I$-bundle over a compact surface.
In particular, $M$ is either a surface 
bundle over the circle with fiber $F$, or
the union of two twisted $I$-bundles glued 
along their common boundary $F$.
\endproclaim

{\it Proof}: Because \caf\ is Reebless, the 
surface $F$ is $\pi_1$-injective [No],
and so lifts to a collection of planes in 
\tim . Their image $C$ in the space
of leaves of
\tcaf\ is a discrete set of points in 
$\bold R$, since it is closed, and any
sequence of distinct points in $C$ with 
a limit point can be used to
find a sequence of points
in $F$ limiting on $F$ in the transverse 
direction, contradicting the
compactness of $F$. The
complementary regions of $F$ in $M$ lift 
to the complementary regions of
the lifts of $F$ in \tim ; in the 
space of leaves they
correspond to the intervals between successive 
points of $C$. Each is bounded
by two points of $C$, and so every component 
\tix\ of the inverse image of a component 
$X$ of $M|F$ has
boundary equal to two lifts of $F$. Because 
\tim\ is simply-connected, as are the
$\partial$-components of \tix , \tix\ is 
simply-connected, and so \tix\ is
the universal cover of $X$.

Because $F$ $\pi_1$-injects into $M$, it 
$\pi_1$-injects into $M|F$, and
hence into $X$. The index of $\pi_1(F)$ 
in $\pi_1(X)$ is equal to the
number of connected components of the 
inverse image of $F$ in the
universal cover of $X$. To see this, 
choose a basepoint \xo\ for $X$ lying
in $F$, and suppose that \gam\ is a 
loop based at \xo\ which is not in the image
of $\pi_1(X)$. Then the lift \tgam\ of 
\gam\ to \tix\ must have endpoints on distinct
lifts of $F$, for otherwise the endpoints 
can be joined by an arc $\widetilde{\alpha}$ in the
lift of $F$, whose projection, since 
\tix\ is simply connected
(so \tgam$*\widetilde{\alpha}$ is null-homotopic)
is a null-homotopic loop
$\gamma * \alpha$ in $\pi_1(X)$. This implies that [\gam] = 
[$\bar{\alpha}$]$\in \pi_1(F)$, a
contradiction. Choosing representatives
from each coset of $\pi_1(F)$ in $\pi_1(X)$, 
and lifting each to arcs with initial
points a fixed lift \txo\ of \xo , we 
find that their terminal points must therefore
lie on distinct lifts of $F$.

But since \tix\ has only two boundary 
components, this means that
$\pi_1(F)$ has index at most two in $\pi_1(X)$. Since $X$ is
irreducible (because $F$ is incompressible 
and $M$ is irreducible), [He, Theorem 10.6]
implies that $X$ is an $I$-bundle over a closed surface.
The resulting description of $M$ follows.
\qed

The foliation by fibers of a bundle over the
circle is always taut. In the other case, when 
$F$ separates, we understand [Br2] the structure of the
foliation \caf\ on each of the two $I$-bundles,
since their boundaries are leaves. If the 
induced foliations can be made
transverse to the $I$-fibers of each bundle, 
then by taking a pair of $I$-fibers, one from each bundle,
and deforming
them so that they share endpoints on $F$, we 
can obtain a loop transverse to the leaves of \caf , so
\caf\ is taut. If the induced foliations cannot 
be made transverse to the $I$-fibers, then $F$ is a torus,
and (after possibly passing to a finite cover) 
the foliation contains a pair of parallel
tori with a Reeb annulus in between. It is then 
straighforward to see that the resulting lifted
foliation \tcaf\ cannot be $\bold R$-covered,
since this torus$\times I$ will lift to $R^2\times I$
whose induced foliation has space of leaves 
$\bold R$ {\it together with} two points 
(the two boundary components)
that are {\it both} the limit of the positive 
(say) ray of the line. In particular, the space
of leaves of \tcaf\ would not be Hausdorff. Therefore:

\proclaim{Corollary 2} A Reebless, $\bold R$-covered 
foliation is taut. \qed\endproclaim

We now turn our attention to $\bold R$-covered 
foliations with Reeb components. Such foliations do exist, for
example, the foliation of $S^2\times S^1$ as a 
pair of Reeb components glued along their
boundaries; the lift to $S^2\times R$ consists 
of a pair of solid cylinders, each having space
of leaves a closed half-line. Gluing the solid 
cylinders together results in gluing
the two half-lines together, giving space of 
leaves $\bold R$. We show, however, that, in some sense, this is
the only such example. Recall that the 
Poincar\'e associate $P(M)$ of $M$ consists of the connected
sum of the non-simply-connected components 
of the prime decomposition of $M$, i.e.,
$M$ = $P(M) \#$(a counterexample to the Poincar\'e Conjecture). 

\proclaim{Lemma 3} If \caf\ is an $\bold R$-covered 
foliation of the orientable 3-manifold $M$, which
has a Reeb component, then $P(M)$ = $S^2\times S^1$. \endproclaim

{\it Proof}: The core loop \gam\ of the Reeb 
component must have infinite order in the fundamental group of
$M$, otherwise the Reeb component lifts to a 
Reeb component of \tcaf ; but since the interior of a Reeb
component has space of leaves $S^1$, this 
would imply that $S^1$ embeds in $\bold R$, a contradiction.

The Reeb solid torus therefore lifts to a 
family of infinite solid cylinders in \tim ,
foliated by planes. The induced foliation
of each closed solid cylinder has space of 
leaves a closed half-line properly embedded in the
space of leaves of \tcaf . Each such half-line 
is disjoint from the others; but since $\bold R$ has only
two ends, this implies that the Reeb component 
has at most two lifts to \tim . This means that the
inverse image of the core loop \gam\ of the Reeb 
solid torus, in the universal cover \tim , consists 
of at most two lines,
and so the (infinite) cyclic group generated 
by \gam\ has index at most 2 in $\pi_1(M)$.
Because $M$ is orientable, it's fundamental 
group is torsion-free, and so by [He, Theorem 10.7],
$\pi_1(M)$ is free, hence isomorphic to 
$\Bbb Z$, and so [He, Exercise 5.3] $P(M)$ is an $S^2$-bundle
over $S^1$. Since $M$ is orientable, this 
gives the conclusion. \qed

Note that the space of leaves in the universal 
cover does not change by passing to finite covers
(there is only one universal cover), and so 
we can lose the orientability hypothesis by weakening
the conclusion slightly. Putting the lemmas together, we get:

\proclaim{Corollary 4} If a 3-manifold $M$ 
admits an $\bold R$-covered foliation \caf , then either
\caf\ is taut or $P(M)$ is double-covered 
by $S^2\times S^1$ . \qed\endproclaim

\heading{\S 2 \\ Taut but not $\bold R$-covered}\endheading

In [BNR, Theorem D], the authors exhibit a 
family of 3-manifolds which admit $C^{(0)}$-foliations
with no compact leaves, but no $C^{(2)}$-foliations 
without compact leaves. Each of the examples
is obtained from two copies $M_1$,$M_2$ of 
(a once-punctured torus)$\times S^1$, glued together
along their boundary
tori by a homeomorphism $A$. What we will 
show is that for essentially the same choices of $A$,
the resulting manifolds admit taut foliations, 
but no $\bold R$-covered ones.

$A$ is determined by its induced isomorphism on
first homology ${\Bbb Z}^2$ of
the boundary torus, and so we will think of 
it as a 2$\times$2 integer matrix
with determinant $\pm 1$. We choose as basis
for the homology of each torus the pair 
(*$\times S^1$,$\partial F\times$*), where $F$ denotes the
once-punctured torus. (Technically, we 
should orient these curves, but because 
all of the conditions
we will encounter will be symmetric with 
respect to sign, the orientations will make no difference,
and so we won't bother.)

Each $M_i$ is a Seifert-fibered space 
(fibered by *$\times S^1$); the manifold $M_A$
resulting from gluing via $A$ is a
Seifert-fibered space iff $A$ glues fiber 
(1,0) to fiber (1,0), i.e., $A$ is upper triangular.
We will assume that this is not the case. 

Let $T$ denote the incompressible torus \del $M_1$ = \del $M_2$
in $M_A$. By [EHN], any horizontal foliation 
of $M_i$, i.e., a foliation everywhere transverse to
the Seifert fibering of $M_i$, must meet 
\del$M_i$ in a foliation with slope in the interval
($-$1,1). [Note that this disagrees with 
the statement in [BNR], where the result was quoted
incorrectly.] 

If \caf\ is an $\bold R$-covered, hence Reebless, 
foliation of $M$, then by [BR], we can
isotope \caf\ so that either it is transverse to $T$,
and the restrictions \cafi\ of \caf\ to 
$M_i$, $i=1,2$, have no Reeb or half-Reeb components,
or \caf\ contains a cylindrical component, 
and therefore a compact (toral) leaf. In the second
case, the torus leaf must hit the torus $T$, 
and is split into a collection of non-\del -parallel
annuli; these (essential) annuli must be 
vertical in the Seifert-fibering of each $M_i$, since
$F\times S^1$ contains no horizontal annuli. 
But this implies that the gluing map $A$ glues
fiber to fiber, a contradiction. Therefore, 
we may assume that \caf\ restricts to Reebless
foliations on each of the manifolds $M_i$. 
Note that [BR] requires that we allow a finite amount of
splitting along leaves to reach this conclusion; 
but since a splitting of an $\bold R$-covered
foliation is still $\bold R$-covered (it amounts, 
in the space of leaves, to replacing points with
closed intervals), this will not affect our argument.

By [Br3], each of the induced foliations \cafi\ 
of $M_i$ has either a vertical or horizontal
sublamination. Every horizontal lamination in 
$M_i$ can be extended to a foliation transverse to the
fibers of $M_i$, and so meets \del $M_i$ in 
curves whose slope lies in $(-1,1)$. If \cafi\
has a vertical sublamination, it either meets 
\del $M_i$ in curves of slope $\infty$ (i.e.,
in curves homologous to (1,0)) or is disjoint 
from the boundary.

It is this last possibility, a vertical 
sublamination disjoint from $T$, that we wish to require,
and so we will now impose conditions on the 
gluing map $A$ to rule out the other possibilities.
If {\it both} \cafu\ and \caft\ have either a 
horizontal sublamination or a vertical
sublamination meeting $T$, then for both $M_1$ 
and $M_2$, the induced foliations meet $T$
in curves with slope in $(-1,1)\cup\{\infty\}$. 
Therefore, the gluing map

\ssk

\ctln{$A$ = $\pmatrix a&b\cr c&d\cr\endpmatrix$ , 
which as a function on slopes,
is A($x$)=$\displaystyle{{ax+b}\over{cx+d}}$}

\ssk

\noindent must have $A((-1,1)\cup\{\infty\})
\cap((-1,1)\cup\{\infty\})\neq\emptyset$. But since

\ssk

\ctln{$\displaystyle A(\infty) = {{a}\over{c}}$, 
$\displaystyle A({{-d}\over{c}})=\infty$,
$\displaystyle A({{1}\over{1}}) = {{a+b}\over{c+d}}$}

\ssk

\noindent and $A$ is increasing if $ad-bc=1$ 
(take the derivative!), we can force \hfill

\noindent $A((-1,1)\cup\{\infty\})$ to be disjoint 
from $(-1,1)\cup\{\infty\}$ by setting

\ssk

\ctln{$|a| >|c|$, $|d| > |c|$, and $\displaystyle 
{{a+b}\over{c+d}} < -1$}

\ssk

For example, we may choose $a = -n$, $b = -nm-1$, 
$c = 1$, and $d = m$, with $n,m\geq 2$, so that

\ctln{$\displaystyle {{a+b}\over{c+d}}$ = 
$\displaystyle -n - {{1}\over{m+1}}$ $<$ $-1$}

\ssk

As the figure below shows, the conditions $-a > c$, 
$d > c$ and $ad-bc = 1$ are also sufficient; what is
needed, essentially, is that neither the graph of 
$A$ nor either of its asymptotes pass through the
square $[-1,1]\times[-1,1]$.

\ssk

\input epsf.tex

\leavevmode

\epsfxsize=3in
\centerline{{\epsfbox{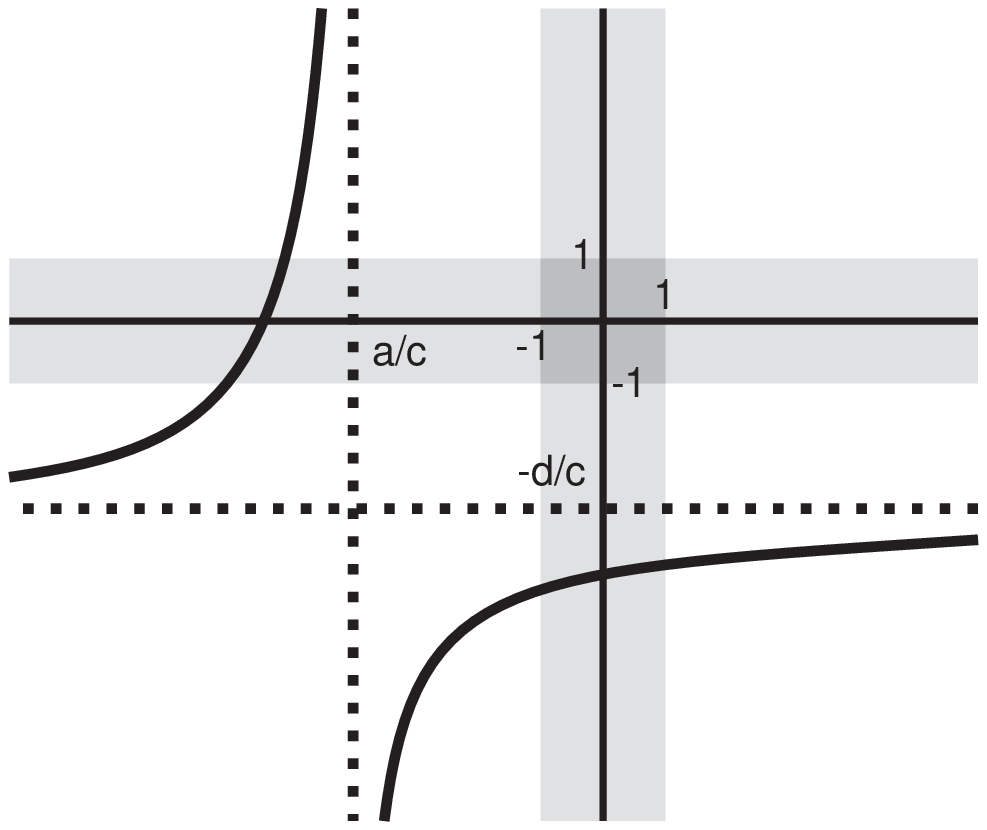}}}

\ctln{Figure 1}

\msk

For such a gluing map $A$ and resulting manifold 
$M_A$, either \cafu\ or \caft\ (without loss
of generality, \cafu ) must
contain a vertical sublamination \cl\ disjoint 
from $T$ = \del$M_1$. \cl\ is the saturation, by circle
fibers, of a 1-dimensional lamination \lam\ in the 
punctured torus $F$. This lamination \lam\ cannot
contain a closed loop, since then \cafu, and 
therefore \caf, would contain a torus leaf $L$ missing $T$.
Lemma 1 would then imply that $M|L$ is an $I$-bundle, 
a contradiction, since it contains $M_2$. \lam\
is therefore a lamination by lines. By Euler 
characteristic considerations, the complementary
regions of \lam, thought as in a torus, are products, 
and so the complentary region of \lam\ in $F$
which meets \del $F$ is topologically a (\del-parallel) 
annulus, with a pair of points removed from the
`inner' boundary. Therefore, the component $N$ of 
$M_A|$\cl\ which contains $T$ is homeomorphic to $M_2$
with a pair of parallel loops removed from \del $M_2$ = $T$ .

Now we will assume that \caf\ is $\bold R$-covered, and argue 
as in the proof of Lemma 1, to arrive at a
contradiction. \cl\ lifts to a lamination in $R^3$ by planes, whose
image in the space of leaves $\bold R$ of \tcaf\ is a closed 
set. The two boundary leaves of \cl\ in \del $N$
are both annuli; they are the complements of the two parallel 
loops in \del$M_2$ in the above description of
$M_A|$\cl . A lift of this complementary region to \tim\ has 
boundary consisting of lifts of the two
annuli. Since the lift is a closed set in \tim\, its image in 
the space of leaves $\bold R$ is a connected, closed
set, and therefore an interval. This implies that the lift of $M_A|$\cl\ 
has (at most) two boundary components,
implying that the inverse image of each of the annulus 
leaves is a single (planar) leaf of \tcaf . This implies, as
in the proof of Lemma 1, that the image in $\pi_1(M_A|$\cl) = $\pi_1(M_2)$ 
of the fundamental group of each
annulus has index at most 1, and so $\pi_1(M_2)$ = 
$F_2\times{\Bbb Z}$ is cyclic, a contradication. Therefore,
\caf\ is not $\bold R$-covered. This implies:

\proclaim{Theorem 5} With gluing map $A$ given as above, $M_A$ admits 
no $\bold R$-covered foliations. \qed\endproclaim

On the other hand, every manifold $M_A$ built out of the pieces we 
have used admits taut foliations and,
in fact, foliations with no compact leaves. We simply choose a vertical 
lamination with no compact leaves
in each of the Seifert-fibered
pieces $M_i$, missing the gluing torus $T$. The complement of this lamination 
is homeomorphic to $T\times I$,
with a pair of parallel loops removed from each of the boundary components. 
Treating this as a sutured manifold, and
thinking of this as
$(S^1\times I)\times S^1$, we can foliate it, transverse to the sutures,
by parallel annuli. Then, as in [Ga1], we can spin 
the annular leaves near the sutures, to extend our vertical laminations to a
foliation of $M_A$ with no compact leaves. This give us:

\proclaim{Corollary 6} There exist graph manifolds admitting taut 
foliations, but no $\bold R$-covered foliations. \qed
\endproclaim

\heading{\S 3 \\ $\bold R$-covered finite covers}\endheading

Work of Luecke and Wu [LW] implies that (nearly) every connected graph manifold is
{\it finitely covered} by a graph manifold that admits an $\bold R$-covered foliation. 
In particular, for
any graph manifold $M$ whose Seifert fibered pieces all have base surfaces having
negative (orbifold) Euler characteristic, they find
a finite cover $M^\prime$ (which is also a graph manifold) 
admitting a foliation \caf\ transverse to the
circle fibers of each Seifert fibered piece of $M^\prime$,
and which restricts on each piece to a fibration over the 
circle. Note that this implies that every leaf
of \caf\ meets every torus which splits $M^\prime$ into Seifert-fibered pieces. 

Even more, every leaf of the lift, to the universal cover 
of $M$, of \caf\ meets every lift $P_1,P_2$
of the splitting tori.
This can be verified by induction on the number of lifts of 
the tori that we must pass through to get from a
lift we know the leaf hits, to our chosen target lift. The 
initial step follows by picking a path \tgam\ between
two `adjacent' lifts
$P_1$ and $P_2$, whose interior misses every lift of the 
splitting tori, and projecting down to $M$; this gives a
path \gam\ in a single Seifert fibered
piece of $M^\prime$. This path can be made piecewise vertical 
(in fibers) and horizontal (in leaves of \caf),
missing, without loss of generality, the multiple fibers of 
$M^\prime$ (just do this locally, in a foliation
chart for \caf; the Seifert fibering can be used as the 
vertical direction for the chart).
Each vertical piece can then be dragged to the boundary 
tori, since the saturation by fibers of an edgemost
horizontal piece of \gam\ is a (singular) annulus with 
induced foliation by horizontal line segments; see Figure 2.
[This is where the fact
that \caf\ is everywhere transverse to the fibers is 
really used.] The end result of this
process is a loop \gampr, homotopic rel endpoints to 
\gam, which consists of two paths each lying in a circle
fiber in the boundary tori, with a single path in 
a leaf of \caf\ lying in between. This lifts to a
path homotopic rel endpoints
to \tgam, consisting of paths in the two lifted tori, 
and a path in some lifted leaf. This middle path
demonstrates that {\it some} lifted leaf $L$ hits 
both $P_1$ and $P_2$.  

\ssk

\leavevmode

\epsfxsize=1.5in
\centerline{{\epsfbox{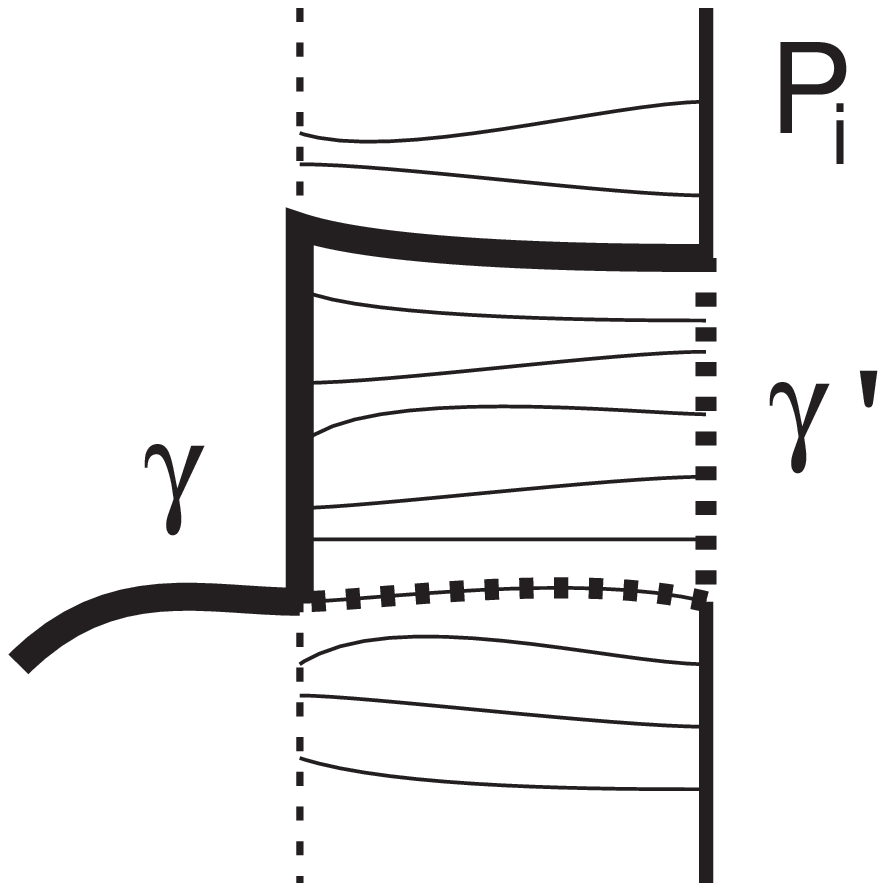}}}

\ctln{Figure 2}

\msk

By choosing a point where any other lifted leaf $L^\prime$ hits a lift $P$ of a
splitting torus, and joining it by a path
$\alpha$ to a point
where $L$ hits $P_1$ or $P_2$, we can apply the same 
straightening procedure as above (see Figure 3), to show that
$\alpha$ is homotopic rel endpoints to paths, one of 
which lies in a lifted fiber and then lies totally in $L^\prime$,
and the other of which lies totally in $L$, and 
then in a lifted fiber. This in particular implies that
$L^\prime$ also hits $P_1$ and $P_2$. Therefore, every 
lifted leaf hits both $P_1$ and $P_2$, as desired.

\ssk

\leavevmode

\epsfxsize=1.5in
\centerline{{\epsfbox{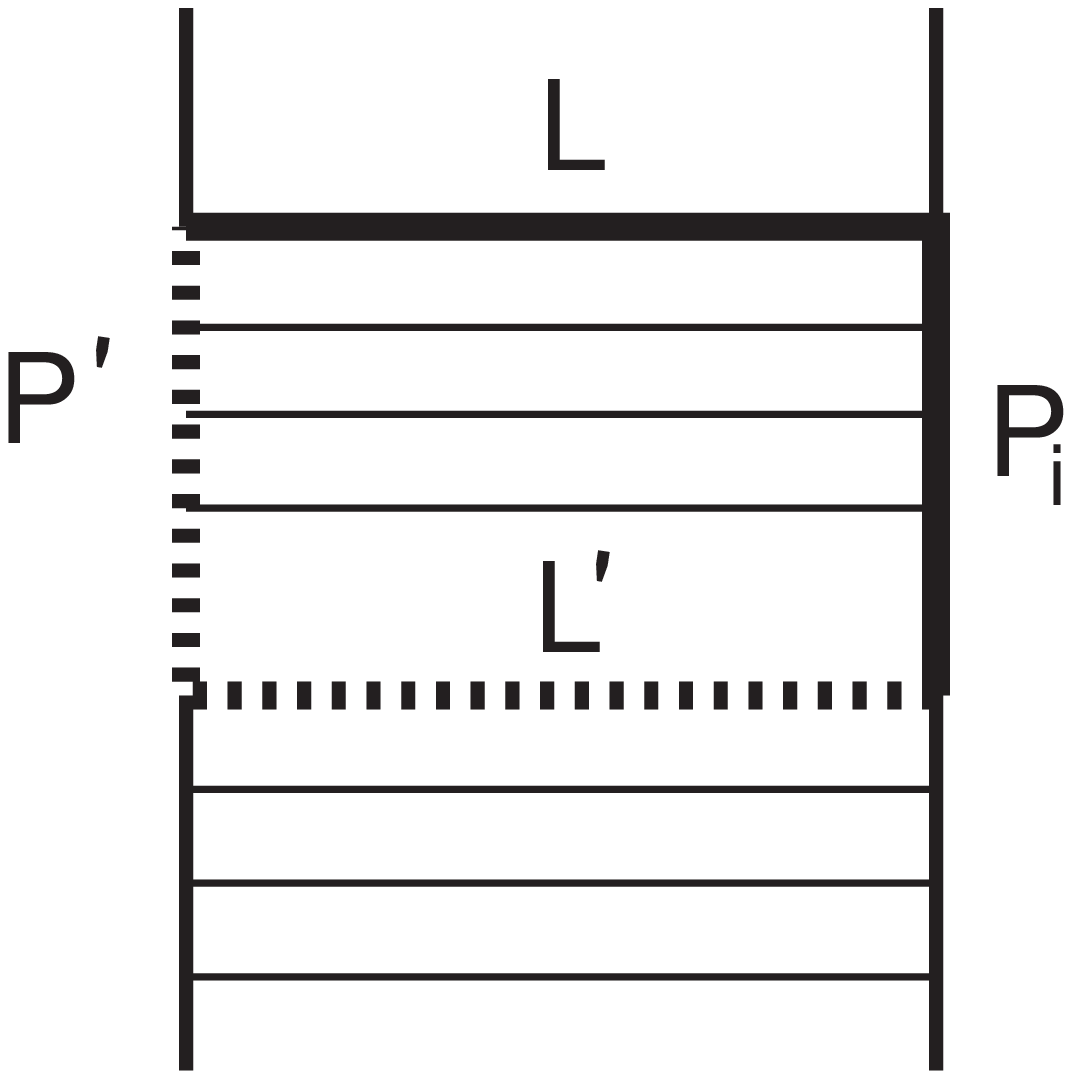}}}

\ctln{Figure 3}

\msk

The inductive step is nearly identical; assuming our 
two leaves $L_1$ and $L_2$ both hit
$P_1, \ldots , P_{n-1}$ , and $P_n$ can be reached 
from $P_{n-1}$ without passing through any other
lift of a splitting torus, the above argument implies 
that two leaves, in the lift of the
relevant Seifert fibered piece, and contained in $L_1$ 
and $L_2$, hit both $L_{n-1}$ and $L_n$, implying
that $L_1$ and $L_2$ also both hit $P_n$.

But this in turn implies that the lifted foliation \tcaf\ has space 
of leaves $\bold R$. This is because the foliation induced by \tcaf\
on any lift $P$ of a splitting torus is a foliation transverse to (either 
of the) foliations by lifts of
circle fibers, and so has space of leaves $\bold R$, which
can be identified with one of the lifts of a circle fiber. [This is 
probably most easily seen in stages:
first pass to a cylindrical cover of the torus, for which the circle fibers lift
homeomorpically. The induced foliation from \caf\ is by lines 
transverse to this fibering, and so has
space of leaves one of the circle fibers. The universal cover 
$P$ is a cyclic covering of this, whose
induced foliation has space of leaves the universal cover of 
the circle fiber.] The argument above implies that
every leaf of \tcaf\ hits $P$ at least once. But no leaf of 
\tcaf\ can hit a lift of a circle fiber {\it more}
than once; by standard arguments, using transverse orientability, 
a path in the leaf joining two such points
could be used to build a (null-homotopic) loop transverse 
to \tcaf, contradicting tautness of \caf, via
Novikov's Theorem [No]. We therefore
have a one-to-one correspondence between the leaves of 
\tcaf\ and (any!) lift of a circle fiber in any
of the Seifert-fibered pieces, giving our conclusion:

\proclaim{Proposition 7} Any foliation of a graph manifold 
$M$, which restricts to a foliation transverse
to the fibers of every Seifert-fibered piece of $M$, is 
$\bold R$-covered. \qed\endproclaim

Combining this with the result of Luecke and Wu, we obtain:

\proclaim{Corollary 8} Every graph manifold, whose 
Seifert-fibered pieces all have hyperbolic base orbifold,
is finitely covered by a manifold admitting an 
$\bold R$-covered foliation. \qed\endproclaim

Combining the proposition with our main resiult, we obtain:

\proclaim{Corollary 9} There exist graph manifolds, 
admitting no $\bold R$-covered foliations, which
are finitely covered by manifolds admitting 
$\bold R$-covered foliations. \qed\endproclaim

\ssk

\heading{\S 4 \\ Concluding remarks}\endheading

Being finitely covered by a manifold admitting an 
$\bold R$-covered foliation is nearly as good as 
having an $\bold R$-covered foliation yourself. Any 
property that could be verified in the presence
of an $\bold R$-covered foliation, which remains 
`virtually' true (e.g., virtually Haken, or having
residually finite fundamental group), would then 
be true of the original manifold. It would then be
of interest to know:

\proclaim{Question 1} Does every 3-manifold with 
universal cover $R^3$ have a finite cover admitting
an $\bold R$-covered foliation?\endproclaim

Or, even stronger:

\proclaim{Question 2} Does every irreducible 
3-manifold with infinite fundamental group
have a finite cover admitting an $\bold R$-covered 
foliation?\endproclaim

Weaker, but still interesting:

\proclaim{Question 3} Does every tautly foliated 
3-manifold have a finite cover admitting an 
$\bold R$-covered foliation?\endproclaim

The first two questions could be broken down into Question 3 and

\proclaim{Question 4} Does every 3-manifold 
(in the appropriate class) have a 
tautly foliated finite cover?\endproclaim

Note that showing that every irreducible 3-manifold 
with infinite fundamental group has a tautly foliated
finite cover would settle the

\proclaim{Conjecture} Every irreducible 3-manifold 
with infinite fundamental group has 
universal cover $R^3$.\endproclaim

Questions 1 and 2 can be thought of as weaker versions 
of the (still unanswered) question, due to Thurston, 
of whether or not every hyperbolic 3-manifold is finitely 
covered by a bundle over the circle; the
foliation by bundle fibers is $\bold R$-covered. Gabai 
[Ga2] has noted that there are
Seifert-fibered spaces for which the answer to Thurston's 
question is `No', although an observation of 
Luecke and Wu [LW]
implies that, via the results [EHN], the answer to our 
Question 1 is `Yes', for 
Seifert-fibered spaces, since [Br4] a transverse 
foliation of a Seifert-fibered space
is $\bold R$-covered. [Note that the arguments of 
the previous section can be modified to give 
a different proof of this;
look at how lifted leaves meet lifts of a single 
regular fiber, instead of lifts of the splitting tori.] 
Question 4, with its conclusion replaced by `have 
a taut foliation', has as answer `No';
examples were first found among 
Seifert-fibered spaces [Br1],[Cl]; there are no known 
examples among hyperbolic manifolds.

\proclaim{Question 5} Do there exist hyperbolic 
3-manifolds admitting no taut foliations?\endproclaim

Finally, the result we have established here for 
graph manifolds is still unknown for hyperbolic
3-manifolds:

\proclaim{Question 6} Do there exist hyperbolic 
3-manifolds which admit taut foliations, but no 
$\bold R$-covered foliations?\endproclaim

\ssk

\enddocument